\documentclass[12pt]{amsart}
\usepackage{graphicx}
\usepackage{graphics}
\usepackage{amsmath}
\usepackage{amscd}
\usepackage{latexsym}

\begin{document}

\textwidth 6.2in
\textheight 7.6in
\evensidemargin .75in
\oddsidemargin.75in

\newtheorem{Thm}{Theorem}
\newtheorem{Lem}[Thm]{Lemma}
\newtheorem{Cor}[Thm]{Corollary}
\newtheorem{Prop}[Thm]{Proposition}
\newtheorem{Rm}{Remark}

\def\a{{\mathbb a}}
\def\C{{\mathbb C}}
\def\A{{\mathbb A}}
\def\B{{\mathbb B}}
\def\D{{\mathbb D}}
\def\E{{\mathbb E}}
\def\R{{\mathbb R}}
\def\P{{\mathbb P}}
\def\S{{\mathbb S}}
\def\Z{{\mathbb Z}}
\def\O{{\mathbb O}}
\def\H{{\mathbb H}}
\def\V{{\mathbb V}}
\def\Q{{\mathbb Q}}
\def\Cn{${\mathcal C}_n$}
\def\CM{\mathcal M}
\def\CG{\mathcal G}
\def\CH{\mathcal H}
\def\CT{\mathcal T}
\def\CF{\mathcal F}
\def\CA{\mathcal A}
\def\CB{\mathcal B}
\def\CD{\mathcal D}
\def\CP{\mathcal P}
\def\CS{\mathcal S}
\def\CZ{\mathcal Z}
\def\CE{\mathcal E}
\def\CL{\mathcal L}
\def\CV{\mathcal V}
\def\CW{\mathcal W}
\def\IC{\mathbb C}
\def\IF{\mathbb F}
\def\IK{\mathcal K}
\def\IL{\mathcal L}
\def\IP{\bf P}
\def\IR{\mathbb R}
\def\IZ{\mathbb Z}

\title{Cappell-Shaneson homotopy spheres are standard}
\author{Selman Akbulut}
\address{Department  of Mathematics, Michigan State University,  MI, 48824}
\email{akbulut@math.msu.edu }
\date{\today}
\begin{abstract} 
We show that an infinite sequence of homotopy $4$-spheres constructed by  Cappell-Shaneson are  all diffeomorphic to $S^4$. This generalizes previous results of Akbulut-Kirby and Gompf.
\end{abstract}

\date{}
\maketitle

\setcounter{section}{-1}

\vspace{-.15in}

\section{Introduction}

Thirty three years ago in \cite{cs}  Cappell and Shaneson defined a sequence of homotopy spheres $\Sigma_{m} $, $m\in {\bf Z}$, as the $2$-fold covers of homotopy ${\bf RP}^4$'s, they constructed (which are known to be exotic when $m=0$ and $m=4$). They asked whether $\Sigma_{m}$  are  $S^4$ or exotic copies of $S^4$. $\Sigma_{m}$  is obtained first by taking the mapping torus of the punctured $3$-torus $T^3_{0}$ with the diffeomorphism induced by the following matrix
\[ A_{m}=
\left(
\begin{array}{ccc}
0  &  1 & 0  \\
 0 &1   & 1  \\
1  &  0 &   m+1
\end{array}
\right)
\]

\noindent  and then by gluing it to a $S^2\times B^2$ with the nontrivial diffeomorphism of $S^2\times S^1$ along their common  boundaries.  In \cite{ak1} it was shown that $\Sigma_{0}$ is obtained from $S^4$ by a ``Gluck construction'' (i.e. by removing a tubular neighborhood of a knotted $S^2$ in $S^4$ then regluing it by the nontrivial  diffeomorphism of $S^2\times S^1$). In \cite{ak1} it was mistakenly claimed that $\Sigma_{0}$ is $S^4$, since at the time we overlooked checking if the gluing diffeomorphism of $S^2\times S^1$ is trivial or not (it turned out  it was in fact nontrivial; this was pointed out in \cite{ar}). Then it took about six years to cancel all the $3$- handles of an handlebody of $\Sigma_{0}$, by turning  it upside down  \cite{ak2}, which resulted a very symmetric  handlebody picture of $\Sigma_{0}$  in Figure 28 of \cite{ak2}, which is  equivalent to Figure 1 ($m=0$ case) below. Then in \cite{g1},  this handlebody of $\Sigma_{0}$ was shown to be diffeomorphic to $S^4$. In \cite{g2}, it was shown that the handlebody of Figure 1 similarly describes $\Sigma_{m}$, for
the cases $m\neq 0$ (see also the discussion in \cite{fgmw}). In this paper we will show the rest of all homotopy spheres $\Sigma_{m}$ are standard (i.e. $m\neq 0$ case).

{\Thm $\Sigma_{m}$ is diffeomorphic to $S^4$, for each $m\in {\bf Z}$.}

\section{The proof}

 \begin{figure}[ht]  \begin{center}  
\includegraphics[width=.60\textwidth]{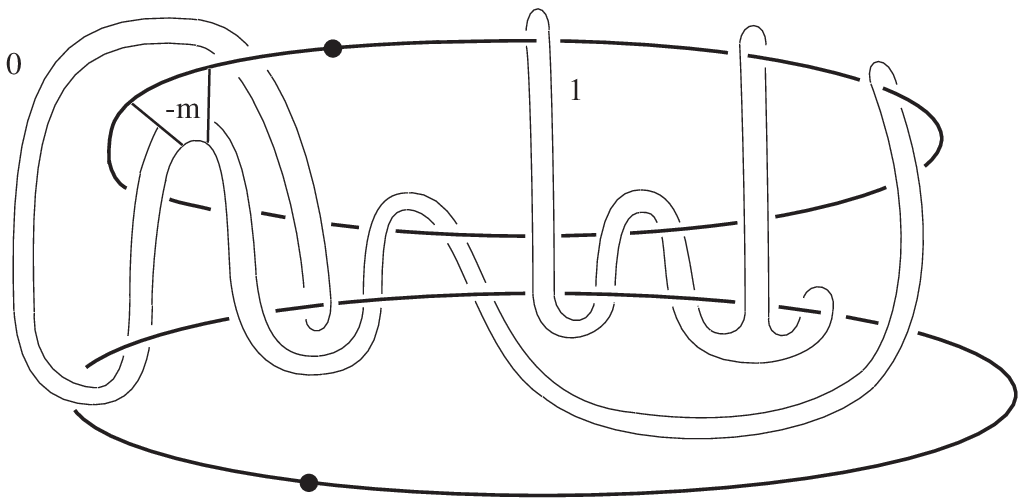}   
\caption{$\Sigma_{m}$} 
\end{center}
\end{figure} 

We first describe a specific  diffeomorphism identifying the boundaries $\partial \Sigma_{m} \approx S^3$. We do this by first surgering the interior of $\Sigma_{m}$ by replacing the two copies of $S^1\times B^3$ with $B^2\times S^2$ (surgery), i.e. we replace the two dotted circles with zero-framed circles.
We then isotope the ``long'' zero framed circle of Figure 1 to the small circle of Figure 2, then surger this zero framed handle (changing the corresponding $S^2\times B^2$ with $B^3\times S^1$). This gives the second picture of Figure 2. Next, we blow down the  $2$-handle corresponding  to the $+1$ framed unknot in the picture (\cite{k}) and get the last picture of Figure 2, which is just the $2$-handle corresponding  to the $-1$ framed unknot in the picture, i.e. it is the punctured $\bar{{\bf CP}}^2$ with boundary $S^3$. 
\vspace{.05in}

Observe that the diffeomorphism $\partial \Sigma_{m} \approx S^3$ shows that the $\alpha$ and the $\beta$ circles on the boundary (Figure 3) are isotopic to each other; each are just $1$ framed unknots in $S^3$. Hence by attaching a $-1$ framed $2$-handle to either $\alpha$, or to $\beta$, we obtain $S^1\times S^2$, and then we can cancel it immediately with a $3$-handle, i.e.  we have diffeomorphisms of the handlebodies (the second and the third handlebodies have $3$-handles):
$$\Sigma_{m} \approx \Sigma_{m} + \alpha^{-1} \approx   \Sigma_{m} + \beta^{-1}$$

Now if we attach $-1$ framed $2$-handle to $\Sigma_{m} $ along $\beta$, by sliding the other $2$-handle going through it, as shown in Figure 4, we see that it becomes just $\Sigma_{m-1}$ with a $2$-handle attached to $\alpha$ with $-1$ framing, i.e.
$$\Sigma_{m} + \beta^{-1} \approx   \Sigma_{m-1} + \alpha^{-1}$$
Hence we have $\Sigma_{m} \approx  \Sigma_{m-1} ...\approx\Sigma_{0}\approx S^4$  \qed.

{\Rm  The Cappell-Shaneson examples come from the 
self diffeomorphisms of $T^{3}_{0}$ induced by a more general family of matrices $A$ up to obvious 
equivalences \cite{cs}, there are finitely many such $A$ for each trace (and only one for each trace between $-4$ and $9$  \cite{ar}). So it was natural to consider this representative family $\Sigma_{m}$ induced by the matrices $A_{m}$. Presumably  there are some more matrices to consider,  but historically authors have been focusing on this   sequence $A_{m}$  (\cite{g1}, \cite{g2}, \cite{fgmw}). To adapt the proof here to other matrices, one has to first construct the handlebody pictures corresponding to Figure 1, which we haven't attempted here. But technique here is not specific to Cappell-Shaneson problem, it is about constructing some hard to see diffeomorphisms between $4$-dimensional handlebodies. In general, beyond guessing,  a useful way to locate such pairs of  $\{ \alpha, \beta \}$ curves is to turn the handlebodies upside down, a technique often used in other  similar problems e.g.  \cite{a1}, \cite{a2}, \cite{a3}.}

\begin{figure}
\includegraphics[width=1.00\textwidth]{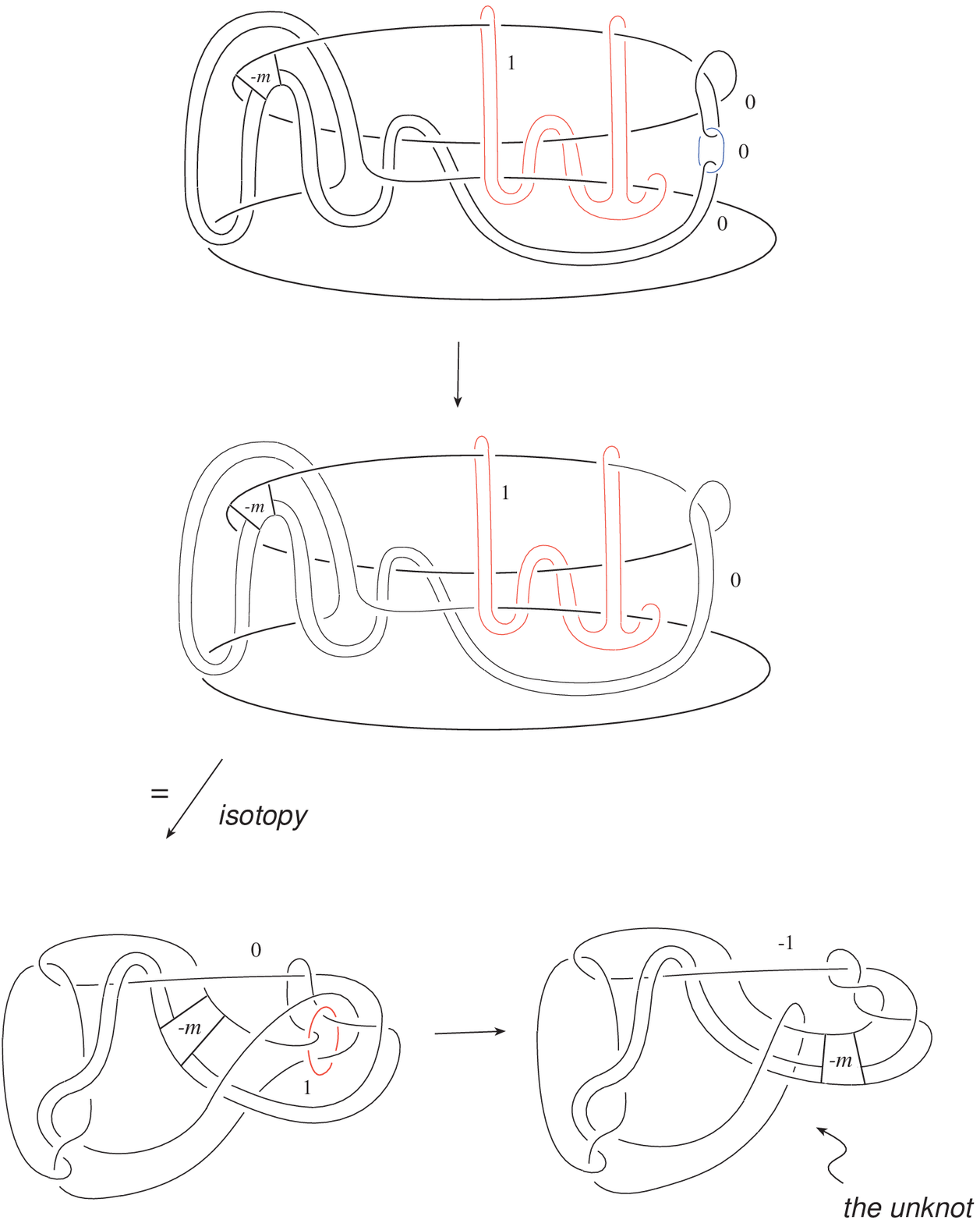}
\caption{}
\end{figure}

\begin{figure}
\includegraphics[width=.80\textwidth]{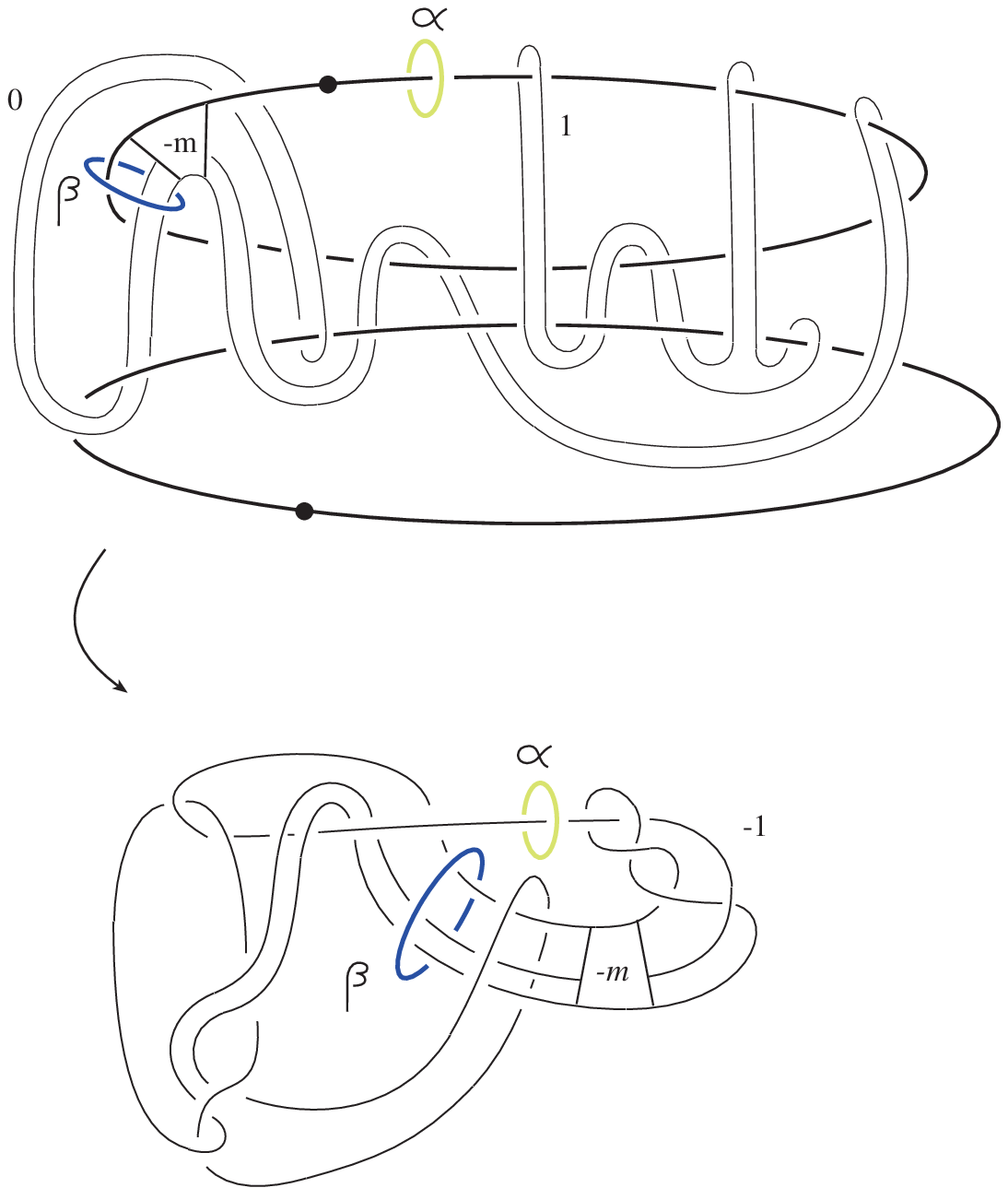}
\caption{}
\end{figure}

\begin{figure}
\includegraphics[width=.70\textwidth]{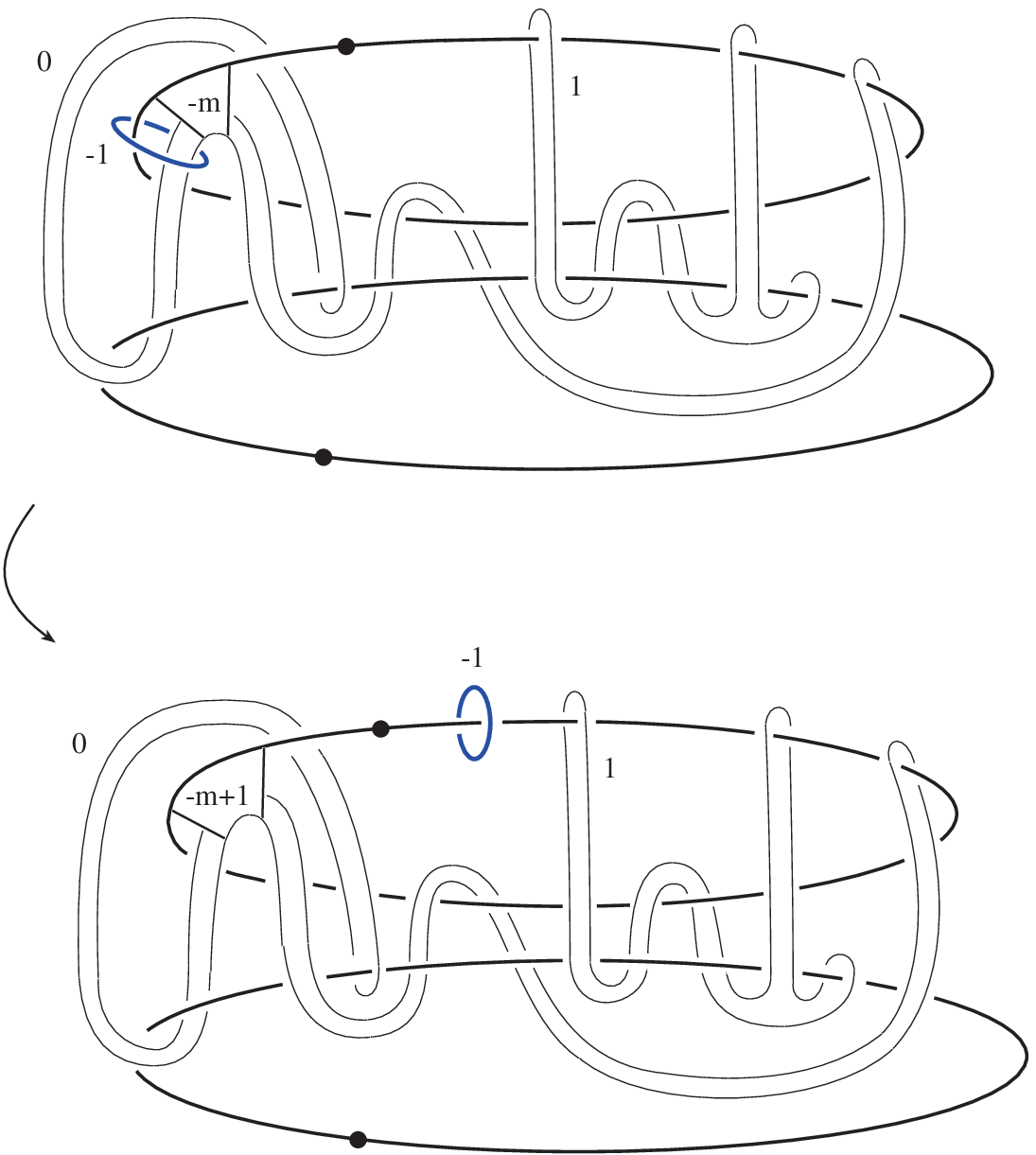}
\caption{}
\end{figure}

\end{document}